# "CÓMO PENSAR COMO SHERLOCK HOLMES". UN MÉTODO PARA ESTUDIAR MATEMÁTICAS

Martha Yoko Takane Imay
Instituto de Matemáticas, Universidad Nacional Autónoma de México
Grupo Mujer y Ciencia, UNAM
yoko.takane@gmail.com

*"El negarle a una niña el placer de acercarse a las Matemáticas, es exactamente una ablación."*
*Martha Yoko Takane Imay*

**Abstract**. *Girls can take many paths to become scientists. But they undoubtedly include the following steps*
*Play, play and play*
*Observe*
*Ask (themselves)*
*In this paper I will talk about some of my experiences teaching girls and boys, teenagers and my students in the science school at UNAM, as well as High school teachers in my course **"How to think like Sherlock Holmes".***

**Resumen**. *Hay muchas maneras de llegar a ser una científica. Pero indudablemente tienen que pasar en algún momento por los siguientes pasos:*

*Jugar, jugar, jugar*
*Observar y*
*Preguntar(se)*

*En este artículo daré algunas de mis experiencias dando clases a niñas y niños, a jovencitas, a mis estudiantes de carreras de la Facultad de Ciencias UNAM y a docentes de Bachillerato en mi curso **"Cómo pensar como Sherlock Holmes"**.*

**Palabras clave:** despertar de la vocación científica, niñas, pensamiento matemático, pensamiento lógico, género, problemas con las Matemáticas, Sherlock Holmes.

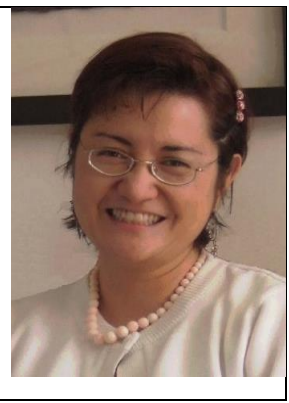

**Algo sobre la autora**: Yoko-niña se la pasaba mirando las nubes, buscando formas e inventando cuentos con ellas. Mirando la luna "redescubrió" el porqué de sus fases y le sigue maravillando ver cómo parece que la sigue. A sus quince años y por el gusto que ya desde entonces tenía por las Matemáticas, una maestra de Biología, amiga de la familia, le regaló los libros "El Hombre que calculaba", "Aritmética recreativa" y "Álgebra recreativa" [4,7,8]. Recuerda también con cariño a sus maestras de Matemáticas y de Química de tercero de secundaria que la introdujeron en el pensamiento científico, donde "redescubrió" el Método REDOX y el pan con mantequilla y azúcar. También recuerda al ejemplo y héroe de toda su vida, su papá, neuro-oftalmólogo de profesión. Ahora es una apasionada del Álgebra y sus aplicaciones.

**Este artículo está organizado de la siguiente manera**: Para entender la razón de este artículo y nuestra lucha para que cambie la manera de cómo se están enseñando, en todos los niveles escolares, las Matemáticas en México, diremos brevemente:

**Qué estudian las Matemáticas y porqué se pueden enseñar y motivar de manera natural.**
Luego daremos un ejemplo de un pedazo de un libro de Sherlock Holmes para ver **cómo se pueden deducir hechos y tener la certeza de lo que se afirma**. En el *Apéndice*, daremos una introducción a la Lógica Proposicional, cuyos ejemplos más ilustres son los Silogismos Aristotélicos y la manera en que Holmes resuelve sus casos. En la *sección 2,* daremos algunas sugerencias de cómo podemos ayudar como mamás, papás, profesoras y profesores a niñas y jovencitas a despertar su vocación científica. Y en general, para toda aquella persona, independientemente de la edad, que quiera gozar de la Ciencia y sus descubrimientos.
Después daremos una lista de ligas interesantes y referencias que creemos interesantes.

**SECCIÓN 1. ¿QUÉ ESTUDIAN LAS MATEMÁTICAS Y CÓMO LO HACEN?**
*Las Matemáticas estudian Todo, pero por medio de abstracciones.*
(La **abstracción** es una operación mental destinada a aislar conceptualmente una propiedad o función concreta de un objeto, y pensar qué es, ignorando otras propiedades del objeto en cuestión. Wikipedia).

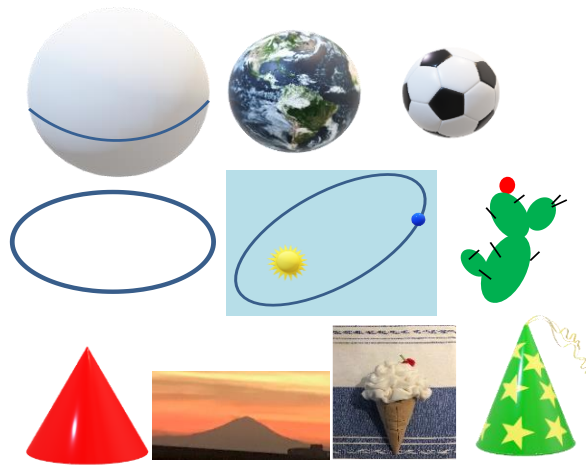

Entre algunos ejemplos de abstracciones se encuentran las siguientes: la abstracción de un balón de futbol y del planeta Tierra es la **esfera**. Estudiamos las órbitas de los planetas y a los nopales, por medio de su abstracción, la **elipse**. Y el **cono** es la abstracción de un cono de helado, de los volcanes y los gorritos de fiesta.
El estudiar objetos naturales por medio de sus abstracciones permite a toda persona que quiera saber propiedades de sus objetos de estudio, por ejemplo, los virus, qué propiedades pueden tener éstos, si por ejemplo son esféricos, así, las Matemáticas inmediatamente dirán qué posibles propiedades, por ejemplo geométricas, pueden tener. Además, estudiar por medio de abstracciones ha dado lugar a inventos importantes, distintas visiones y métodos de estudio en varias ciencias; además de la creación de nuevas áreas científicas como el estudio de redes (sociales y neuronales), entre otras.

Las ideas de la Lógica Matemática que se consideraban en su época desquiciantes, dieron origen a la Computación. Recordemos a Alan Turing (1912-1954) el padre de la programación moderna y a Joan Clarke (1917-1996) criptóloga que con sus investigaciones matemáticas dieron origen a las ciencias de la computación.

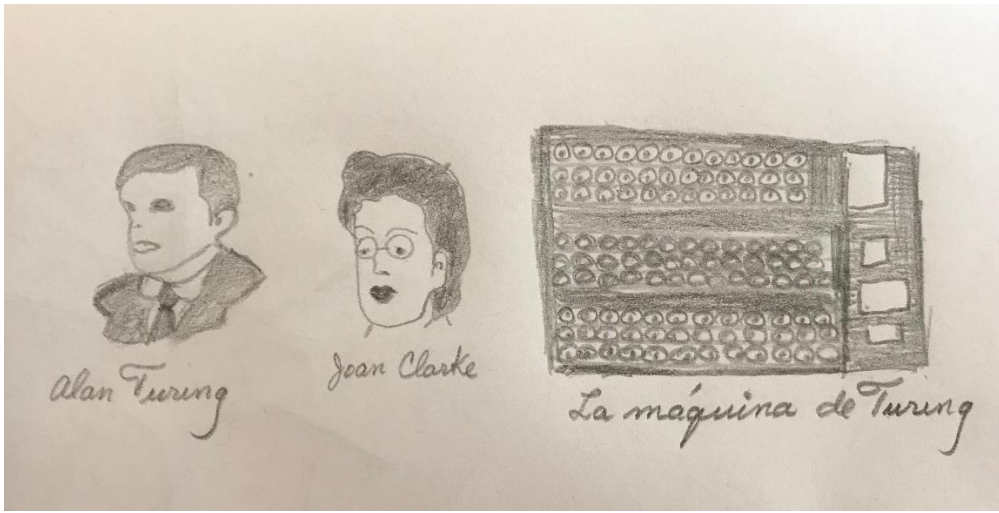

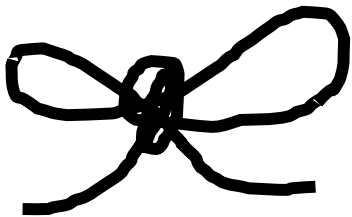 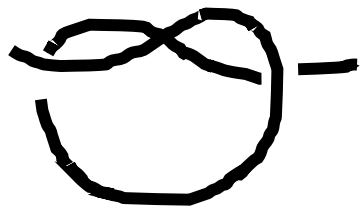 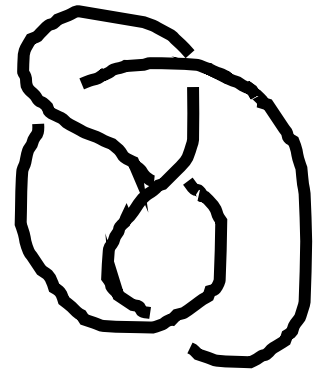

NUDOS

**LOS NUDOS**. Las cadenas de ADN (Ácido Desoxirribonucleico) que forman las proteínas no están, en la naturaleza, como siempre nos las dibujan: como escalera larga-larga... Estas cadenas se anudan y más aún, también los distintos nudos dan propiedades a las moléculas formadas. Por lo que, saber qué nudo es exactamente el que tiene una molécula es muy importante y mejor aún, si no tenemos que romper el nudo para saber cuál es. Esto también lo estudian las Matemáticas y lo resolvieron.

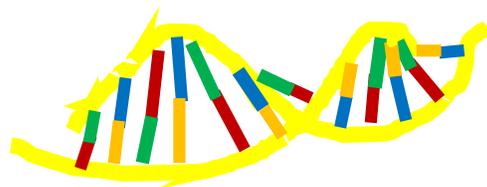

Cadenas de ADN

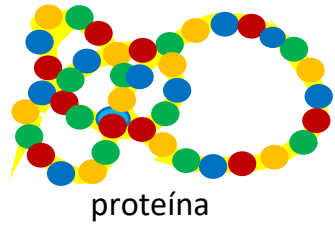

proteína

**LAS REDES**: Carreteras, Rutas de aviones, Redes Neuronales (cerebro, sistema nervioso), Redes sociales (Tus Amig@s, Contactos de Email, Facebook, Google).

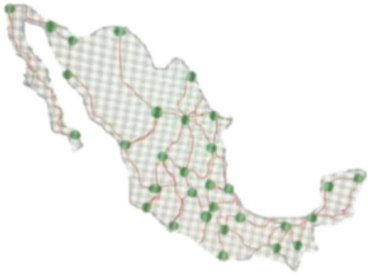

Por ejemplo, se utiliza en redes de transporte para ver qué ciudades están mejor conectadas por medio de distintos transportes y ver en dónde es conveniente y más barato poner un aeropuerto. También se pueden usar en casos de emergencia para poner cercos sanitarios. En el estudio del funcionamiento de las redes neuronales, entre muchas otras aplicaciones.

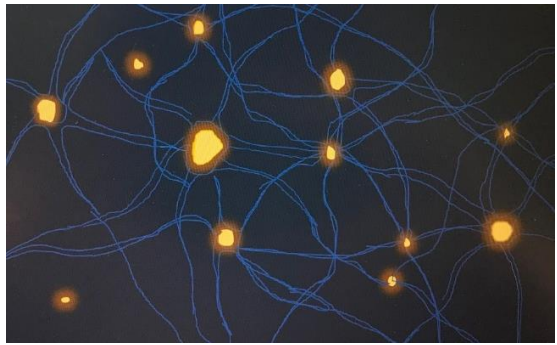

Red Neuronal

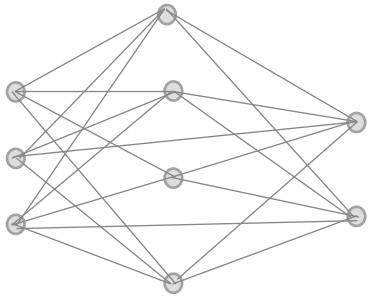

La abstracción de cualquier red es la gráfica

**CÓMO PENSAR COMO SHERLOCK HOLMES.**

Veamos un fragmento del libro [2] "The return of Sherlock Holmes: The Adventure of the Dancing Men" de Sir Arthur Conan Doyle, para ver cómo deduce Holmes sus casos:

"-Así pues, Watson -dijo de repente-, no tiene usted la menor intención de invertir fondos en esas acciones sudafricanas. Di un respingo de asombro. Por muy acostumbrado que estuviera a las extrañas facultades de Holmes, aquella intrusión repentina en mis pensamientos más íntimos resultaba totalmente inexplicable. -¿Y cómo puede usted saberlo? - le pregunté.
-Vamos, Watson, confiese que he logrado asombrarle.
-Claro que sí.
-Debería hacérselo firmar.
-¿Por qué?
-Porque dentro de cinco minutos me dirá usted que era de una sencillez absurda.
-Estoy seguro de que no diré nada semejante.
-Vamos a ver mi querido Watson -empezó a disertar con aspecto de profesor que imparte su clase, después de meter la probeta en su sitio. -No es tan difícil establecer toda una serie de deducciones, cada una de las cuales se apoya en la anterior y es muy sencilla de por sí. Cuando se ha llevado a cabo esa progresión, si se enuncia el resultado desprovisto de todas las deducciones intermedias, y si se presenta al auditorio solamente el punto de partida y la conclusión, no dejará de causar un efecto pasmoso, aun cuando sea algo falso. Ahora bien, no ha sido tan difícil examinando el surco que hay entre su índice y su pulgar, asegurarse de que no tenía usted la intención de arriesgar su pequeño capital en esas minas de oro.
-No veo qué relación tiene una cosa con otra.
-Es probable que no, pero puedo demostrarle a usted, muy rápidamente, una relación muy estrecha. Aquí tenemos los eslabones que faltan en esa cadena tan sencilla:
1) tenía usted manchas de tiza entre el índice y el pulgar izquierdos cuando regresó anoche del club: 2) usted se pone tiza ahí cuando juega al billar, para que resbale bien el taco, 3) sólo juega usted al billar con Thurston, 4) me dijo usted hace cuatro semanas que Thurston tenía una opción sobre alguna propiedad sudafricana que expiraría al cabo de un mes, y que deseaba compartir con usted, 5) tiene usted el talonario de cheques guardado en mi cajón y no me ha pedido la llave, 6) así pues, no tiene usted intención de arriesgar su dinero en esas condiciones.
-¡Qué absurda sencillez! - exclamé.
-Exactamente - respondió Holmes, un poco ofendido- Todos los problemas resultan infantiles una vez explicados ...".

**AHORA DAREMOS UN EJEMPLO SENCILLO DE CÓMO TRABAJAMOS: Desde un problema cotidiano a la abstracción.**

## "EL PUESTO DE AGUAS FRESCAS"

Dos amigas, *A* y *B*, quieren poner un puesto en la Feria para vender aguas frescas. Mientras *B* pone el puesto, *A* va a comprar dos envases, con el encargo de que uno sea de medio litro y el otro de un litro.
**Un tiempo después:** *A* llega con dos vasos y dos helados. *A* estaba muy contenta pues los vasos estaban de barata y hasta le alcanzó a comprar un helado para cada una. El problema es que un vaso era de 3 dl. y otro de 6 dl. *B* la mira con desconcierto: "Pero queríamos vender litros y medios litros."

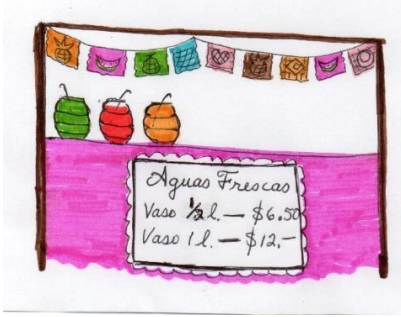

**PROBLEMA A RESOLVER:** ¿Qué cantidades *exactas* de aguas frescas pueden vender, ahora? Necesitan hacer el cartel, ¿qué cantidades ponen?
**Aquí empieza el reto matemático.**
Empecemos **experimentando (jugando)** con ambos vasos y hagamos observaciones**:**

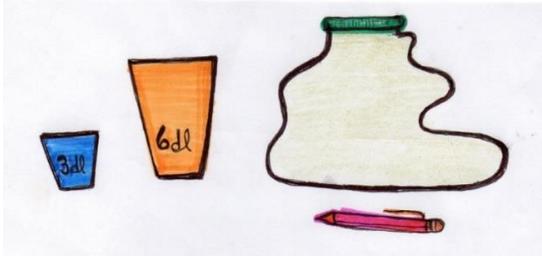

Para esto buscan cualquier envase, suficientemente grande, y un plumón para marcar las cantidades de agua que se pueden hacer con ambos vasos.

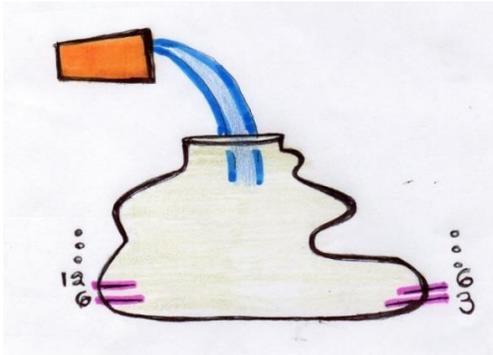

Se dan cuenta que se pueden hacer: 3dl., 6dl. ,9dl., 12dl. y todas las cantidades que son múltiplos de 3. También se dan cuenta de que pueden hacer 6dl., 12dl., 18dl. y cualquier cantidad que sea múltiplo de 6. Y como el 6 es múltiplo de 3, es decir, 6=3X2, se dan cuenta que con el vaso de 3dl. pueden hacer todas las cantidades que se pueden hacer con el de 6dl.

POR LO QUE AHORA, EL CARTEL TENDRÍA QUE DECIR:

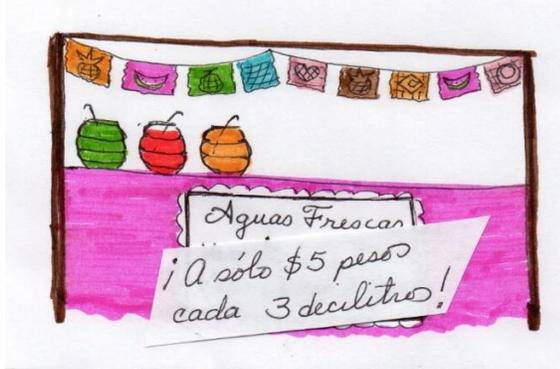

TAMBIÉN, al estar jugando, se dan cuenta de lo siguiente:

¡SE PUEDE RESTAR!
¡SI PODEMOS HACER EL 1 (UNO), PODEMOS HACER TODOS LOS NÚMEROS! (pero en nuestro caso, sólo podemos hacer múltiplos de 3).
¡DE HECHO, NO IMPORTA SI TENEMOS DECILITROS, LITROS, MILILITROS, GALONES, SE TENDRÍAN LOS MISMOS NÚMEROS EN NUESTRAS CUENTAS!

HAGAMOS OTRO EXPERIMENTO: **¿qué cantidades podremos hacer si tenemos un vaso de 3dl. y 11dl.?**
COMO ANTES, PERO AHORA YA QUITANDO "dl.", PODEMOS HACER:
3,6,9,…. Todos los múltiplos de 3
11,22,33,…. Todos los múltiplos de 11.

Y TAMBIÉN EL UNO!!

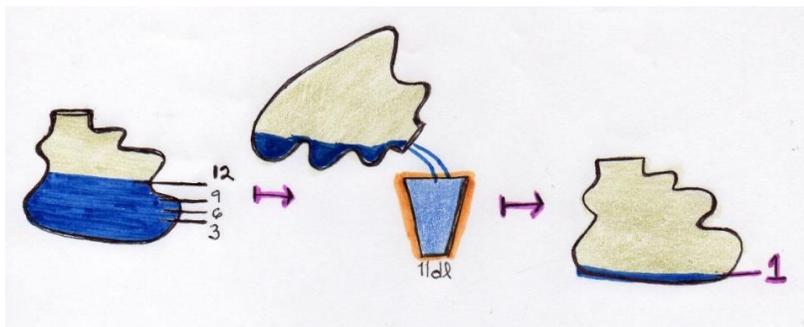

Es decir, si ponemos, en el envase que conseguimos, 4 medidas del vaso de 3dl. y con eso llenamos el vaso de 11dl., lo que nos resta en el envase es exactamente 1dl. A saber, **1**=12-11=(**3**X4)-**11.**

**AHORA, A LA ABSTRACCIÓN**: Tomemos cualesquiera dos cantidades, digamos **n** y **m** que como observamos anteriormente, pueden ser de cualquier unidad de medida. Entonces podemos construir, todos los múltiplos de **n**, todos los múltiplos de **m** y más aún, todas las cantidades que sean de la forma **an+bm** donde **a** y **b** son números enteros, que como podemos restar, pueden ser negativos, pero como queremos que sean cantidades de agua, tenemos que construir todos los **an+bm** que sean mayores que 0, cero.

Más aún, las Matemáticas nos dicen, que estos números, **an+bm, son LOS MÚLTIPLOS DEL máximo común divisor de n y m.**
**¡Y con eso podemos resolver completamente nuestro problema!**

SIGUIENDO NUESTROS EJEMPLOS:
1. En el caso de 3dl. y 6dl., el máximo común divisor de 3 y 6 es 3. Y como vimos, pudimos hacer con estos vasos, exactamente, es decir ni más ni menos, todos los múltiplos de 3.

Recordemos:
Los divisores del 3 son 1,3.
Los divisores del 6 son 1,3,6.
Entonces los números que dividen al 3 y al 6, sus divisores comunes son 1,3. Por lo tanto, el máximo de ellos es 3.

2. **El Máximo común divisor del 3 y el 11 es 1,** pues

Los divisores del 3 son 1,3
Los divisores del 11 son 1,11
Entonces sólo tienen al 1 como divisor común. Por lo tanto, el máximo común divisor del 3 y el 11 es el 1.

<div align="center">¡LAS MATEMÁTICAS SE PUEDEN USAR MÁS DE LO QUE SE CREE!</div>

**SECCIÓN 2. ¿CÓMO AYUDAR A NUESTRAS HIJAS, NUESTRAS ALUMNAS, JOVENCITAS PARA QUE SE EMPODEREN ANTE LA VIDA Y EN PARTICULAR A ESTUDIAR MATEMÁTICAS O ALGUNA CIENCIA EXACTA? En las referencias agregamos ligas interesantes que podrían ser útiles para este objetivo.**

1.Como dijimos anteriormente, necesitamos que aprendan o lo hagamos junto con ellas a
   Jugar, jugar y jugar
   Observar
   Y Preguntar(se)
2. Tener un diccionario de palabras y frases que NO se deben decir ni en la casa ni en la escuela porque agreden la autoestima de la niña.
3. Tener un diccionario de palabras y frases que SÍ se pueden decir y que ayuden a la autoestima de la niña.
4. Y se vale que digamos, como profesoras, profesores, mamás y papás: "No sé, pero lo buscaremos y aprenderemos junt@s".
5. Crear microambientes sanos.
6. Solidaridad y generosidad de TODOS los niños, jóvenes y hombres adultos.
Por ejemplo, si en los medios de transporte, hay un vagón o carro especialmente para niñas y mujeres, respetarlo. Algunos hombres me han dicho que eso va en contra de su libertad de tránsito.
Necesitamos su solidaridad, para que mientras México sea un lugar inseguro, podamos protegerlas de alguna forma.
También solidaridad con nuestras niñas y jovencitas, de que si un amigo hombre está haciendo algo indebido, sus amigos hombres le digan que no lo haga.
7. Algo que se ve que ha servido, es dar talleres o clases de Matemáticas sólo para niñas o jovencitas, dónde se puedan desarrollar con libertad.

8. Sí hay historias de éxito en México, por ejemplo, las Olimpiadas Mexicanas Femeniles de Matemáticas.

**¡ES BUEN MOMENTO PARA SER UNA CIENTÍFICA, NO les inculquemos miedo, por el miedo que tenemos como madres y padres a que fracasen, NO FRACASARÁN!**

## Apéndice. Introducción a la Lógica Proposicional.

**Conceptos primitivos: Verdadero (*V* ó 1) y Falso (*F* ó 0).** (Que llamaremos "**valores verdad**").

**Definición**. Una **proposición** es un enunciado (oración, frase) que es verdadera o falsa PERO no ambas.

Fijemos ideas: Tomen como verdadero o falso lo que usualmente signifiquen para ustedes.

**Ejemplos y no ejemplos de proposiciones:**
1. "El número 5 es un número impar". Esta oración es verdadera, por lo tanto, es una proposición.
2. "El número 2 es un número impar". Esta oración es falsa, por lo tanto, es una proposición.
3. "El número 9". Esta frase no es una proposición.
4. "¡Vámonos!", no es una proposición.

**Notación.** Denotaremos a las proposiciones con letras mayúsculas *P, Q, R, S, T, …*

**Las reglas del juego (maneras de construir nuevas proposiciones) o AXIOMAS de la lógica proposicional son los siguientes:**

Dadas dos proposiciones *P* y *Q*, podemos formar nuevas proposiciones, más aún, dependiendo de los valores verdad de *P* y *Q*, los axiomas nos dirán qué valor verdad tendrá la nueva proposición construida.

**AXIOMAS de los conectivos**: Sean *P* y *Q* proposiciones. Las siguientes también serán proposiciones ¬*P* (ó "no es cierto que pasa P"), ***P ó Q*** *(pasa P ó pasa Q")*, ***P y Q*** *("pasa P y pasa Q")*, ***P*⟹*Q*** *("Si pasa P, entonces pasa Q")* y ***P*⟺*Q*** *("Pasa P exactamente cuando pasa Q").*

La **Negación**: **"No es cierto que pasa *P*"**, que denotaremos por ¬***P*.**
Si *P* es verdadera entonces ¬ *P* es falsa. Y si *P* es falsa, entonces ¬ *P* es verdadera.
Es decir, su **tabla verdad** (o tabla de valores verdad) es la siguiente:

| P | ¬P |
|---|---|
| V | F |
| F | V |

**Ejemplo.** Sea *P* la proposición *" El número 5 es impar",* que sabemos es verdadera*.*
Su negación ¬*P* es *"No es cierto que el número 5 es impar"* o *"El número 5 no es impar",* que es (y debería de ser) falsa.

La **Disyunción**: ***P ó Q*** que se lee "pasa *P* o pasa *Q*".
Esta proposición es verdadera si alguna de las proposiciones (o ambas) es verdadera.

Para construir la tabla verdad de *P ó Q*, necesitamos tener todas las combinaciones de valores verdad de *P* y *Q*. A saber,

| P | Q | P ó Q |
|---|---|---|
| V | V | **V** |
| V | F | **V** |
| F | V | **V** |
| F | F | **F** |

**Ejemplo**: Sea *P* la proposición "*El número 5 es impar*", que es verdadera y *Q* la proposición "*El número 8 es impar*", que es falsa. Por lo tanto, la proposición *P ó Q* es verdadera. Es decir, la proposición "*El número 5 es impar o el número 8 es impar*" es verdadera.

La **Conjunción**: *P y Q* que se lee "pasa *P* y pasa *Q*".
Esta proposición es verdadera si ambas *P* y *Q* son verdaderas, en los demás casos será falsa.

| P | Q | P y Q |
|---|---|---|
| V | V | **V** |
| V | F | **F** |
| F | V | **F** |
| F | F | **F** |

**Ejemplo**: Sea *P* la proposición "*El número 5 es impar*", que es verdadera y *Q* la proposición "*El número 8 es impar*", que es falsa. Por lo tanto, la proposición *P y Q* es falsa. Es decir, la proposición "*El número 5 es impar y el número 8 es impar*" es falsa.

El **Si... entonces... (o implica o condicional)**: *P$\Rightarrow$Q* que se lee "si pasa *P* entonces pasa *Q*".
Esta proposición es la que se usa para generar conocimiento nuevo "verdadero", por lo que observemos cuidadosamente su tabla verdad, donde el único caso en que el implica es falso es cuando *P* es una proposición verdadera y *Q* es falsa. Es decir, si empezamos con algo verdadero NUNCA implicaremos algo falso. PERO CUIDADO, si empezamos con algo falso, podemos llegar a proposiciones verdaderas o falsas, es decir, *no podemos confiar en este nuevo conocimiento*. Así actúa la realidad, que es lo que tratamos de modelar, y claro, a veces (que no debería de pasar) se manipula la información. Por lo que SIEMPRE debemos de tener la certeza de empezar con proposiciones verdaderas.

| P | Q | P$\Rightarrow$Q |
|---|---|---|
| V | V | **V** |
| V | F | **F** |
| F | V | **V** |
| F | F | **V** |

**Ejemplo de cómo generamos conocimiento**: Sea *P* la proposición *"Está lloviendo en mi colonia"* y *Q* la proposición *"El pasto de mi jardín está mojado"*.

Ahora, por la experiencia o por alguna otra razón, sabemos que la proposición *"Si está lloviendo en mi colonia, entonces el pasto de mi jardín está mojado"* es verdadera, es decir, la proposición $P \Rightarrow Q$ es verdadera. Además, ahorita está lloviendo, es decir, *P* es verdadera. Entonces la única manera en que tengamos que *P* es verdadera y $P \Rightarrow Q$ es verdadera es que *Q* sea verdadera y ya generamos conocimiento*: Tengo la seguridad de que mi jardín está mojado.*

Observemos que mi jardín puede estar mojado por otras razones, por ejemplo, que mi hermana regó el jardín. Pero, $P \Rightarrow Q$ lo que dice es que si pasa *P* entonces, seguro, pasa *Q.*

El **Bicondicional (ó "pasa… exactamente cuando pasa…" ó "si y sólo si")**: $P \Leftrightarrow Q$ que se lee "pasa *P* exactamente cuando pasa *Q*".

Esta proposición es verdadera cuando los valores verdad de *P* y *Q* son idénticos.

| P | Q | $P \Leftrightarrow Q$ |
|---|---|---|
| V | V | **V** |
| V | F | **F** |
| F | V | **F** |
| F | F | **V** |

Ejemplo: Sean *P* la proposición "*a+b=a+c*" y Q la proposición "*b=c*", donde *a,b* y *c* son números enteros.
Entonces $P \Leftrightarrow Q$ es la proposición "*a+b=a+c* si y sólo si *b=c*".

**Definiciones**. Sea *P* una proposición. Entonces
   a) *P* es una **tautología** si su tabla verdad tiene s lo valores Verdaderos.
   b) *P* es una **contradicción** si todos sus valores verdad son Falsos.

**Ejemplos:** La proposición *P ó ¬P* es una tautología y la proposición *P y ¬P* es una contradicción, como veremos a continuación*.*

| P | ¬P | *P ó ¬P* | *P y ¬P* |
|---|---|---|---|
| V | F | **V** | **F** |
| F | V | **V** | **F** |

**Axioma de los Cuantificadores**. Sea *P(x)* una proposición que involucra a *x*. Por ejemplo, sea *P(x)* la proposición "*x* es un número par".

Entonces podemos hacer nuevas proposiciones de las siguientes maneras:

**Cuantificador Universal**. La siguiente es una proposición "Para todo *x*, pasa *P(x)*".
   Esta proposición es verdadera si para cada *x* la proposición *P(x)* es verdadera.

**Cuantificador Existencial.** La siguiente es una proposición "Existe *x* tal que pasa *P*(*x*)" o "Hay (un) *x* tal que cumple *P*(*x*)" o "Hay (al menos un) *x* tal que (pasa) *P(x)*".

Esta proposición es verdadera si exhibimos una *x* tal que la proposición *P*(*x*) es verdadera.
Negaciones de los cuantificadores:
¬( Existe *x* tal que pasa *P*(*x*)) es la proposición Para todo *x*, no es cierto que pasa *P*(*x*).
¬( Para toda *x* pasa *P*(*x*)) es la proposición Existe *x* tal que no es cierto que pasa *P(x)*.

**Definición**. Sean *P* y *Q* dos proposiciones. Diremos que *P* y *Q* son (proposiciones) **equivalentes** si *P*⇔*Q* es una tautología.
**Importante**: Si *P* y *Q* son proposiciones equivalentes, entonces en cualquier lugar donde aparezca *P* podemos intercambiarla por *Q* y viceversa.

**Observación**: Algunas tautologías son tan importantes que tienen nombres, como las siguientes:
1. *Modus ponens*: **[(*P*⟹*Q*) y *P*]** ⟹*Q*. Y se lee: Si pasa *P* entonces pasa *Q* y (además) pasa *P*, entonces pasa *Q*.
2. *Tollendo ponens*: **[(*P* ó *Q*) y (¬*P*)]** ⟹*Q*.
3. *Tollendo tollens*: **[(*P*⟹*Q*) y (¬*Q*)]** ⟹ (¬*P*).
4. *Contrapuesta*: **[*P*⟹*Q*]** ⟹ **[(¬*Q*) ⟹ (¬*P*)]**.
5. *Silogismo hipotético*: **[(*P*⟹*Q*) y (*Q*⟹*R*)]** ⟹ **[*P*⟹*R*]**.
6. *Dilema constructivo*: **[(*P*⟹*Q*) y (*R*⟹*S*) y (*P* ó *R*)]** ⟹ **[*Q* ó *S*]**.
7. *Dilema destructivo*: **[(*P*⟹*Q*) y (*R*⟹*S*) y [(¬*Q*) ó (¬*S*)]]** ⟹ **[(¬*P*) ó (¬*R*)]**
8. *Exportación*: **[*P*⟹ (*Q*⟹*R*)]⇔[(*P* y *Q*) ⟹*R*]**

**LEMA**: Sean *P* y *Q* proposiciones. Supongamos que son equivalentes, entonces *P* y *Q* tienen exactamente la misma tabla verdad.

**SILOGISMOS ARISTOTÉLICOS**. A continuación daremos algunos de estos silogismos, es decir, proposiciones que son tautologías de la forma *(P y Q)⟹R.* Estos silogismos tienen nombre, para recordar cuáles son.

**Bárbara**: Si todo M es B y todo A es M entonces todo A es B.
**Celarent**: Si ningún M es B y todo A es M entonces ningún A es B.
**Darii**: Si todo M es B y algún A es M entonces algún A es B.
**Ferio**: Si ningún M es B y algún A es M entonces algún A no es B.
**Cesare**: Si ningún B es M y todo A es M entonces ningún A es B.
**Camestres**: Si todo B es M y ningún A es M entonces ningún A es B.
**Festino**: Si ningún B es M y algún A es M entonces algún A no es B.
**Baroco**: Si todo B es M y algún A no es M entonces algún A no es B.

**Darapti**: Si todo M es B y todo M es A entonces algún A es B.
**Felapton**: Si ningún M es B y todo M es A entonces algún A no es B.

**REFERENCIAS Y LIGAS INTERESANTES PARA AYUDAR AL DESPERTAR DE LA VOCACIÓN CIENTÍFICA EN NIÑAS Y JOVENCITAS.**

**JUEGOS para todas las edades**

Professor Fizzwizzle
Es un juego donde un profesor tiene que ir de un lugar a otro, pero hay obstáculos. Es muy interesante porque se parece mucho a cómo trabajamos en Matemáticas
Interactive Mathematics Miscellany and Puzzles: http://www.cut-the-knot.org/

**REFERENCIAS Y LIGAS INTERESANTES SOBRE IGUALDAD Y EQUIDAD DE GÉNERO Y SOBRE LENGUAJE INCLUYENTE.**
https://mujeryciencia.wixsite.com/unam
https://proyectos.matem.unam.mx/igualdad

https://www.matem.unam.mx/
https://www.museodelamujer.org.mx/
https://www.museodelamujer.org.mx/index.php?page=24  (Enlaces interesantes)
http://www.welcomingschools.org/resources/school-tips/transgender-youth-what/trans-how/este-preparado-para-preguntas-e-insultos-sobre-el-genero/
https://twitter.com/juanitolibritos/status/1217084113435471872?s=12

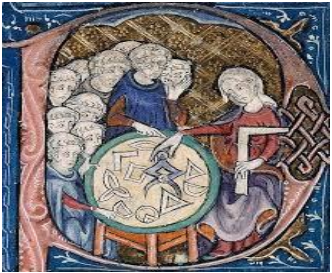
Mujer enseñando geometría (1309-1316)
Fuente: The British Library

**1er. Congreso Internacional "El Despertar de la Vocación Científica en las Niñas"**
**29, 30 y 31 de octubre de 2018**
**Sedes: Instituto de Astronomía e Instituto de Matemáticas, UNAM**
**https://www.matem.unam.mx/fsd/takane**

Para referir este artículo: Memorias del 1er. Congreso Internacional "El Despertar de la Vocación Científica en las Niñas", Universidad Nacional Autónoma de México (UNAM).